\documentclass{article}
\usepackage{graphicx}
\usepackage{amssymb}
\usepackage{a4}
\newtheorem {Theorem}{Theorem}

\newtheorem{Problem}{Problem}

\def\Aut{{\rm Aut}}

\def\ad{{\mathrm{ad}}}

\def\lie{{\mathfrak l}}
\def\la{{\mathfrak a}}

\def\lg{{\mathfrak g}}
\def\lh{{\mathfrak h}}

\def\CC{{\mathbb C}}

\title{Algorithmic Construction of Chevalley Bases}
\author{ K. Magaard 
\\ School of Mathematics,\\ 
University of Birmingham,\\
Birmingham B15 2TT,\\ U.K.\\
and \\
 R. A. Wilson\\
School of Mathematical Sciences,\\ 
Queen Mary, University of London,\\ 
Mile End Road,\\ London E1 4NS,\\ U.K.}

\begin{document}

\maketitle
\begin{abstract}
We present a new algorithm for constructing a Chevalley basis
for any Chevalley Lie algebra over a finite field. 
This is a necessary component for some constructive recognition
algorithms of exceptional quasisimple groups of Lie type.
When applied to a simple Chevalley Lie algebra in characteristic $p\ge 5$,
our algorithm has complexity
involving the 7th power of the Lie rank, which is likely to be
close to best possible.
\end{abstract}
\section{Introduction}

Finding a Chevalley basis for a semisimple
Lie algebra over $\CC$ amounts to diagonalizing
a regular semisimple element: the eigenspaces for non-zero
eigenvectors are just the $1$-dimensional root spaces, and suitable eigenvectors
can be chosen as described by Carter \cite{Carter}. 
Indeed, the same is true for any Chevalley Lie algebra
over any algebraically closed field.
However, over a finite field the problem is much more difficult. 
The probability that a random
regular semisimple element is split is approximately the reciprocal of the order of
the Weyl group, so something better than a random search is required 
if we want a polynomial-time algorithm.

Let us define a {\em toral subalgebra} of a Lie algebra $\lie$ to be any abelian
subalgebra consisting of semisimple elements.
If $\mathfrak t$ is a maximal toral subalgebra which is split, then
its centralizer in $\mathfrak l$ is a \emph{Cartan subalgebra}
$\mathfrak c$, and conversely,
$\mathfrak t$ consists exactly of the semisimple elements in $\mathfrak c$.

\begin{Problem}\label{prob1}
Given a 
split toral subalgebra $\lh_0$ in a Chevalley Lie algebra
$\lie$, find a 
Cartan subalgebra $\lh$ such that $\lh_0 \subset \lh$,
and a Chevalley
basis with respect to $\lh$.
\end{Problem}
\begin{Problem}\label{prob3}
 Given two
Cartan subalgebras $\lh_1, \lh_2$ of a semisimple Lie algebra $\lie$,
find an element $g \in \Aut(\lie)$
such that $\lh_1^g = \lh_2.$
\end{Problem}

Solutions to these problems are a necessary component for some constructive recognition
algorithms of exceptional quasisimple groups of Lie type \cite{KM}.
A polynomial-time Las Vegas algorithm for solving Problem \ref{prob1} 
is given by Ryba \cite{Ryba},
except in characteristic $2$
(where indeed Problem \ref{prob1} has no solution in general), 
and except for $\la_2$ and $\lg_2$ in characteristic $3$.
This algorithm has complexity involving the $11$th power of the Lie rank of the
algebra,
as well as the fourth power of the logarithm of the field order, 
although practical implementations are apparently much faster than this suggests.
He asserts that his algorithm often works in characteristic $2$, but does not attempt a
full analysis in that case.

Another algorithm is given by Cohen and Murray \cite{CM}, with the same exceptions,
with complexity (in the case when the input is an algebra corresponding to a simple
algebraic group) involving the 9th power of the Lie rank. (A noteworthy feature of
their algorithm is that the rate-determining step seems to be checking at each
stage whether they have finished. It is possible therefore that their algorithm can
be improved by a more subtle approach to this particular step.)
They do not discuss the exceptional cases. 

The small characteristic
exceptions are discussed by Cohen and Roozemond \cite{CR}, but
they only consider 
the problem of finding a Chevalley basis once a 
Cartan subalgebra has been found. They do not solve the problem of
finding such a subalgebra in the first place.
(This problem is dealt with by Roozemond in \cite{Rooz}.)
Problem \ref{prob3} amounts to finding a base-change matrix which maps one
Chevalley basis to another, so is easily reduced to Problem \ref{prob1}, 
as will be discussed at the end of Section~\ref{algorithm}.

In this paper we propose a 
simpler algorithm which has 
better complexity than 
the above algorithms in the simple case.
We achieve this
by computing the whole Chevalley basis at once, rather than by first
computing the Cartan subalgebra. 
We build up the Dynkin diagram one
node at a time, making each connected component in full before moving
on to the next. Our main theorem is as follows.

\begin{Theorem}
\label{mainthm}
Let $\lie$ be a Chevalley 
Lie algebra 
over a field of order $q$ and 
characteristic $p\ge 5$. Suppose $\lie$ has Lie rank $l$ and dimension $d$. Then there
is an algorithm to compute a Chevalley basis of $\lie$ in $O(ld^3\log q)$ field
operations.
\end{Theorem}

\section{The main algorithm}
\label{algorithm}
We assume that
the characteristic of the field is at least $5$.
In this case our strategy is to look for
a (long or short root) fundamental $\mathfrak a_1$, and find its Chevalley basis $\{e,f,h\}$.
Then we look for another fundamental $\mathfrak a_1$ which extends it to
a simple rank $2$ algebra (if there is one). Continuing in this way, we build up the
connected component of the
Dynkin diagram one node at a time. Then we
iterate the procedure until all components are dealt with.
 Once all components are completed, we use
the `extraspecial pairs' as described by Carter \cite{Carter} to complete the Chevalley
basis for the corresponding simple Lie subalgebras.
The algorithm in detail is as follows. (Comments on `suitable' choices follow the algorithm.)

\begin{enumerate}
\item Input:  a Chevalley Lie algebra $\lie_0$ over a finite field
of characteristic $p\ge 5$, and a split toral subalgebra
$\lh_0$ (defaulting to zero).
\item Output: a Cartan subalgebra $\lh$ containing $\lh_0$, together with the part 
 of 
 a Chevalley basis 
for $\lie_0$, consisting of the $e_\alpha,f_\alpha, h_\alpha$ for simple roots $\alpha$, and
 a complete weight space decomposition $\mathcal W$ of $\lie_0$.
\item 
Initialise $\lh_1:=0$. Initialise $\lh:=0$.
Initialise $\mathcal D:=\emptyset$.
\item If $\lh_0\ne 0$, compute the weight spaces for $\lh_0$, and set $\mathcal W$ equal to
the set of weight spaces, and pair the weight-spaces for opposite non-zero weights. 

 Else pick a random $x\in \lie_0$ and compute the eigenspaces of
$\ad x$ on $\lie_0$, until there are some non-trivial eigenspaces with non-zero eigenvalues,
and set $\mathcal W$ equal to this set of eigenspaces, paired as before.
Adjoin to $\mathcal W$ the perp of $\mathcal W$, so that $\mathcal W$
spans the whole space.
\item  Until $\mathcal W$ consists of a single subspace which is abelian,
\begin{enumerate}
\item  Using the current $\mathcal W$,
find an $\mathfrak a_1$ subalgebra, as follows:
\begin{enumerate}
\item Until there is a pair of opposite $1$-dimensional members of $\mathcal W$, 
pick a pair of opposite spaces $V^+,V^-\in\mathcal W$ with $\dim V^+$ minimal,
and pick random $y\in V^+$ and $z\in V^-$, and let
 $x=[y,z]\in[V^+,V^-]$, and refine the members of $\mathcal W$
using the eigenspaces
of $\ad x$ and the perp. Recompute the pairing of 
members of $\mathcal W$. 
\item
Pick $e\in V^+$ and $f\in V^-$, 
and set $h=[e,f]$. Then scale $h$ so that $[h,e]=2e$; then scale
$f$ so that $[e,f]=h$. Set $\lh_1:=\langle h\rangle$.
Set $\mathcal D:=\{e,f,h\}$.
\item Compute the eigenspaces of $\ad h$. 
Refine $\mathcal W$ using these eigenspaces.
Label each element of
$\mathcal W$ by the corresponding eigenvalue of $\ad h$. (This label is the first coordinate
of what will become the weight vector.)
\end{enumerate}
\item Repeat until 
a maximal string diagram has been found.
\begin{enumerate}
\item  Repeat: pick a suitable label ($\lh_1$-weight) $w$ 
where the next node of the diagram might live, and
a pair $V^+,V^-$ of opposite spaces in $\mathcal W$, with labels $\pm w$,
and random $x\in[V^+,V^-]$, and compute eigenspaces
of $\ad x$ on  $V^+$ and $V^-$;
 until $\ad x$ has a pair of
$1$-dimensional eigenspaces $\langle e\rangle\subseteq V^+$, 
$\langle f\rangle\subseteq V^-$ 
for non-zero
eigenvalues $\pm\lambda$.
\item  Set $h=[e,f]$. Then scale $h$ so that $[h,e]=2e$; then scale
$f$ so that $[e,f]=h$.
\item Adjoin $h$ to $\lh_1$. 
 Adjoin to $\mathcal D$
the vectors $e,f,h$.
\item Compute the eigenspaces of $\ad h$.
Refine $\mathcal W$ using these eigenspaces.
Append to the label of each element of
$\mathcal W$ the corresponding eigenvalue of $\ad h$.
\end{enumerate}
\item Analyse the string diagram obtained in the previous step, to see whether or not it is equal to the
Dynkin diagram of the current component, using the data and notation from Table~\ref{weightdims}. 
\begin{enumerate}
\item If the diagram has just two nodes, then for 
both end nodes, compute $\dim V_1$, $\dim V_2$ and $\dim V_3$
to determine both what the diagram is and what it 
should be:
the only case where it could be wrong, is when the diagram is $A_2$ but should be $G_2$.
\item Else compute $\dim V_2$ for both end nodes
 and one
interior node. 
\item If one of these is $> 1$ and $\ne 7$,
then the diagram is $B_n$ or $\tilde C_n$, and if $B_n$, is correct; if $\tilde C_n$, 
delete an end node to obtain $C_n$.
\item If one of these is $7$,
the diagram is $F_4$ or $B_4$. Distinguish these by considering the node adjacent to
the short end node. If it is $F_4$, the diagram is correct.
Otherwise, compute $\dim V_1$ for the short end node: if this is $0$, the diagram is correct.
Otherwise, it is a $B_4$ diagram but should be $F_4$.

\item Else all nodes of the diagram are long, and the diagram is $A_n$. 
Compute $\dim V_1$ to determine what the diagram should be.
The possible cases where the diagram should not be $A_n$ are as follows:
$D_{n+1}$ (any $n$); $n=2, G_2$; $n=3, D_k (k\ge 4)$; $n=4, E_8$, $n=5, E_6, E_7$;
$n=7, E_7, E_8$; $n=8, E_8$. 
 \end{enumerate}

\item If the type of the diagram is not 
what we know it should be, 
adjust it 
to be the diagram of the current component,
as follows, correcting $\lh_1$, $\mathcal W$ and $\mathcal D$ as we go:
\begin{enumerate}
\item $A_2$ instead of $G_2$: repeat the above steps until roots of both
lengths are found.
\item $B_4$ instead of $F_4$: remove the long end node, and attach a new (short) node
at the other end.
\item $A_n$ instead of $D_{n+1}$: adjoin a node to the penultimate node.
\item $A_3$ instead of $D_n$: attach a tail to the middle node.
\item $A_m$ instead of $E_n$: attach a tail to a suitable node.
\end{enumerate}
\item 
Write out $\mathcal D$, all of $\mathcal W$ which consists of $1$-spaces
labelled by non-zero weights, together with these labels.

Adjoin $\lh_1$ to $\lh$.

Remove the part of $\mathcal W$ which has been written out, and initialise labels
to $\emptyset$. Initialise $\mathcal D:=\emptyset$.

\item
 If  $\mathcal W=[W]$, pick a random $x\in W$ and compute the eigenspaces of
$\ad x$ on $W$, until there are some non-trivial eigenspaces with non-zero eigenvalues,
and set $\mathcal W$ equal to this set of eigenspaces, paired as before.
(If this fails, then $W$ is probably abelian, so break.)
Adjoin to $\mathcal W$ the perp of $\mathcal W$, so that $\mathcal W$
spans the whole space.
\end{enumerate}
\item 
Now $\lh$ is a subspace of the single element of $\mathcal W$, so adjoin
to $\lh$ a complement. Write out $\lh$.
\end{enumerate}

\begin{table}
\caption{Dimensions of weight spaces for $h$ in a simple Lie algebra}
\label{weightdims}
$$\begin{array}{cc|ccc|}
\mbox{Type}& \mbox{Root} &\dim V_1 & \dim V_2 & \dim V_3\cr
\hline
A_n & & 2(n-1) & 1&\cr
D_n & & 4(n-2) & 1&\cr
E_6&&20 & 1&\cr
E_7& &32 & 1&\cr
E_8& &56 & 1&\cr
B_n & \mbox{short} &0 &2n-1 &\cr
B_n & \mbox{long} & 4n-6& 1&\cr
C_n & \mbox{short}& 4(n-2)& 3&\cr
C_n & \mbox{long} & 2(n-1)& 1&\cr
F_4 & \mbox{short}& 8 & 7&\cr
F_4 & \mbox{long}& 14 & 1 &\cr
G_2 & \mbox{short}&2&1&2\cr
G_2 & \mbox{long}& 4&1&\cr
\hline
\end{array}$$
Note: $V_\lambda$ denotes the eigenspace with eigenvalue $\lambda$. Since
$\dim V_{-\lambda}=\dim V_\lambda$, we omit half the eigenvalues. We assume that
the characteristic of the field is at least $7$: obvious modifications of this table apply
in smaller characteristics.
\end{table}

\paragraph{Comments on the algorithm.}
In Step 5(ii)(a), we construct the string diagram by repeatedly trying to attach a node
to the previous one. 
This means looking in the weight space corresponding to the weight $(0,\ldots,0,1)$
or $(0,\ldots,0,2)$. It is clear from Table~\ref{weightdims} that the former case pertains
except in the case $B_n$ at the first step, if a short root has been found.
When this process terminates, we reverse the string, (and the corresponding orderings of $\mathcal D$ and $\lh_1$, and the labels on elements of $\mathcal W$)
and try again,
using the same weights. When this also terminates,
a maximal string subdiagram has been found.

In Step 5(iv)(e) there are various cases, which we now describe in more detail.
In each case the $A_m$ diagram we have found is a maximal $A_m$
subdiagram of the extended
Dynkin diagram $\tilde E_n$.
If $m=4$ (so $n=8$), we attach a tail of length $4$ to one of the two interior nodes:
a priori, we do not know which, so try both. If $m=5$ and $n=6$, attach a node to
the middle node. If $m=5$ and $n=7$, attach a tail of length $2$ to a node
adjacent to one of the  end nodes: again we do not know in advance which one, so try both.
If $m=7$ and $n=7$, try to attach a node to the middle node, after removing either one
of the end nodes. If $m=7$ and $n=8$, we try the following procedure for each
node adjacent to an end node: first remove the far end node, and then attach a tail
of length $2$. Finally, if $m=n=8$, do the same for the nodes at distance $2$
from each end.

\paragraph{The Chevalley basis.}
At completion of the main algorithm, 
we have obtained a Cartan subalgebra, and a complete set of
root vectors for the fundamental roots and their negatives. We also have 
a set of vectors which are scalar multiples of the other root vectors.
It remains to complete this to a Chevalley basis of the commutator
subalgebra by computing the correct
scalar multiples of these.

We assume that for every abstract
Dynkin diagram, a choice of structure constants has been made
(see Chapter 4 of \cite{Carter}). Then we scale each $e_{\alpha+\beta}$ in turn
to ensure that $[e_{\alpha},e_{\beta}]$ is the appropriate multiple ($0,\pm1,\pm2,\pm3$) of
$e_{\alpha+\beta}$. This requires the characteristic to be at least $5$ in the case of 
a component $G_2$,
and at least $3$ in the cases of a component $B_n,C_n,F_4$. 
In each case, to compute the scalar, it suffices to compute one non-zero coordinate of
$[e_{\alpha},e_{\beta}]$. 
This can be accomplished by computing just one column of $\ad e_{\beta}$ and
applying it to $e_{\alpha}$. 
Once all these scalars have been computed, we have a complete Chevalley basis for
$[\lie_0,\lie_0]$.

\paragraph{Solution to Problem~\ref{prob3}.} 
In the case when $[\lie,\lie]=\lie$ we may use our algorithm
with input $\lh_1$ to produce a Chevalley
basis containing a basis of $\lh_1$, and again with input $\lh_2$.
Then any linear map which takes the first basis to the second,
preserving the labelling of the root system,
will be an automorphism of the algebra mapping $\lh_1$ to $\lh_2$, as required.

\section{Analysis of the algorithm}
We first analyse the algorithm and its complexity in the case when the
input algebra is simple and no partial Cartan subalgebra is given.

Let $l$ be the Lie rank, and $d\sim l^2$ the dimension of the algebra,
and let $q$ be the order of the field.

Computation of $\ad x$ for a random vector $x$ takes $O(d^3)$ field operations.
To compute a pair of eigenspaces for non-zero eigenvalues $\pm\lambda$ 
(which we do not compute), we use \cite{CellerLG}, which takes $O(d^3\log q)$ 
field operations.
Computing $[x,y]$ also takes $O(d^3)$ field operations, for example by
computing $\ad y$ and applying it to $x$.

At the start of the algorithm (Step 4) 
we are looking for an element $x$ such that
$\ad x$ has a pair $\pm \lambda$ of non-zero eigenvalues.
The proportion of such elements is at least a constant, say $1/3$
 (see Corollary 6.3 of \cite{CM}). Hence this step can be accomplished in
$O(d^3\log q)$ field operations.

In the simple case the main loop (Step 5)
will be traversed only once. 
In Step 5(i)(a), (and similarly in Step 5(ii)(a))
the commutator $[y,z]$ is in effect a random matrix of small rank 
in the centralizer of the part of the Cartan subalgebra that we have seen. 
The statistics of this situation are at least as good as the
statistics for a random element.
Thus Step 5(i)(a) takes a constant number of $O(d^3\log q)$ steps.
To justify Step 5(i)(b) we need to show that $e$ and $f$ generate a split $\mathfrak a_1$
subalgebra. This follows from the Jacobi identity for $x,e,f$ and for $x,h,e$.
That is
$$[h,x]=[[e,f],x]=[[e,x],f]+[e,[f,x]]=[\lambda e,f]-[e,\lambda f]=0,$$
and
$$[[h,e],x]=[[h,x],e]+[h,[e,x]]=0+[h,\lambda e]$$
so $[h,e]$ is a $\lambda$-eigenvector of $\ad x$ so is a scalar multiple of $e$.
Hence, from the representation theory of $\mathfrak{sl}_2$
we know  in particular that $\ad h$ is diagonalisable.
Thus Step 5(i)(b) works, and takes a constant number of $O(d^3)$ steps. 
Moreover, the eigenvalues of $\ad h$ lie in $\{0,\pm1,\pm2,\pm3\}$
so its eigenspaces can also be computed in $O(d^3)$ field operations.

Step 5(ii)(a) is done (at most) once for each fundamental root, and the computations
each time are essentially the same as in Step 5(i).
Hence this takes $O(ld^3\log q)$ field operations.
Step 5(ii)(b) consists of at most a constant number of eigenspace
computations for known eigenvalues, so takes $O(d^3)$ operations.
Step 5(ii)(c) is similar to 5(ii)(a), and might be done $O(l)$ times if we were in the
case where we mistook $D_l$ for $A_3$.

Steps 5(v) and 5(vi) are book-keeping and termination so do not take
significant time.

The final step of computing the scalars for each weight space for a non-simple root
takes $O(d^2)$ field operations for each root.
Thus this computation can be done in time $O(d^3)$.

Hence the overall complexity in the simple case is $$O(ld^3\log q)
=O(l^7\log q)=O(d^{3.5}\log q)$$ field operations.
The proof of Theorem~\ref{mainthm} is now complete.

\section{Non-simple algebras}
\label{reductive}
\paragraph{Semisimple case.}
We have designed our algorithm to apply to the semisimple case,
by ensuring that in Step 5(i) we at least halve the dimension every time
we find a new eigenspace. Hence this step needs to be applied at most
$\log d$ times to find an $\la_1$ in the first component. 
Since each application of Step 5(i) or Step 5(ii) reduces the rank by $1$,
the overall complexity becomes
$O(ld^3\log d\log q)$.

\paragraph{Non-trivial centres.} The part of the centre which is generated by
commutators is part of the output of the algorithm. The rest of the centre
plays no role, and we can pick an arbitrary basis for it.

\paragraph{Imperfect algebras.}
In this case, extra non-central toral elements
 appear in the final step of the algorithm. 
However, in general it is not possible to scale these to any
particularly nice form. For example, such an $h$ may act non-trivially on multiple
components, and it is only possible to scale it to act canonically on one component.
If the derived subalgebra has large homogeneous components and large codimension,
this makes the definition of a canonical basis almost impossible.

In certain cases, however, it is possible to extend our algorithm. For example,
if the derived subalgebra is simple, then there is at most one dimension
of non-central torus outside the derived subalgebra, and we can make a canonical
choice of element. For example, we can demand that $[e_i,h]=\delta_{i1}e_i$,
where $e_i$ correspond to the fundamental roots. 

\section{Characteristics $2$ and $3$}
\label{char23}
\paragraph{Characteristic $3$.}
The main problem in small characteristics is that in certain cases the weight spaces are
not $1$-dimensional. There may be additional problems for small fields.
In characteristic $3$
we only encounter problems with multidimensional weight spaces
 in the cases where the Lie algebra
has a component of type $\mathfrak g_2$, or a simply-connected component of
type $\mathfrak a_2$. In both these cases, there are eigenspaces of dimension
$3$. Consider first the case $\mathfrak g_2$. In this case, the short roots
occur in weight spaces of dimension $1$, so these are obtained with high
probability in the same way as above, i.e. by looking for a short root $\mathfrak a_2$.
Then we need to modify the algorithm in Step 5(ii)(b) to test whether this $\mathfrak a_2$
should be a $\mathfrak g_2$: specifically, we compute the image of $\ad e$ for one of the
short roots $e$, and test whether this lies in the $\mathfrak a_2$ algebra.
If it does not, then we deduce that the algebra generated by the $\mathfrak a_2$ and this
image is the full $\mathfrak g_2$, so modify Step 5(ii)(c) accordingly, using \cite{CR}.

The simply-connected $\mathfrak a_2$ case
will only arise at the end, when we have run out of $1$-dimensional eigenspaces,
and only $3$-dimensional eigenspaces remain. For each pair of these, we compute the
algebra they generate, and find a suitable basis using \cite{CR}.
See also \cite{Rooz}. We expect that these modifications will not affect
the overall complexity of our algorithm.

The only other problem in characteristic $3$ is in Step 5(ii)(b), where we cannot distinguish
easily between long and short roots in $F_4$ using Table~\ref{weightdims}.
In this case we may have picked up a $B_4$ root subsystem rather than the whole $F_4$.
In order to detect this, we need to check directly whether all $48$
root vectors lie in the algebra generated by the fundamental root vectors.
If not, then we can correct the fundamental roots in the same way as before.

\paragraph{Characteristic $2$.}
We expect that a combination of our ideas with those of \cite{CR} and \cite{Rooz} will
also produce a more efficient algorithm in characteristic $2$. First we briefly sketch
how this might work in the simple   case $A_n$.

\begin{enumerate}
\item Take random $x$, until we have a $2$-dimensional
 eigenspace of $\ad x$ with nonzero eigenvalue.
Pick $e,f$ at random in this eigenspace until $h=[e,f]\ne 0$.
\item Find an eigenspace $V$ of $\ad h$ with non-zero eigenvalue, and
scale $f$ and $h$ so that the eigenvalue is $1$.
\item  Let $V_e=[V,e]\cap V$ and $V_f=[V,f]\cap V$, and pick $y\in V_e$, $z\in V_f$
until $x:=[y,z]\ne 0$. 
\item $\ad x$ acts on $V_e$ and $V_f$, so
intersect the eigenspaces of $\ad x$ with $V_e$ and $V_f$.
Similarly for $\ad(x+h)$.
This gives us enough $1$-dimensional spaces to define the root spaces
for an $\mathfrak a_2$ subalgebra. Scale the vectors as far as possible.
\item Continue in this way to generate each node of the diagram in turn.
\end{enumerate}

More generally, there is no pairing of weight spaces, and the minimal
dimension eigenspaces which we are aiming for have dimension at most $8$
(see \cite[Table 1]{CR}).
If we modify Step 5(i)(a) by taking $V^+=V^-\in\mathcal W$ then we will
reach such a small-dimensional eigenspace in at most $\log d$ steps.
If this dimension is not $2$ or $4$ then the component is of bounded rank,
and the methods of \cite{CR} suffice.
In the other cases,
 we can analyse the subalgebra generated by this eigenspace
in the same way as in \cite{CR}, or as suggested above in the dimension $2$ case.
We then exend to the whole component by a modified version of
Step 5(ii)(a): we know which eigenspace $V=V^+=V^-$ 
to look in, and
if this has dimension $2$ we proceed as suggested in Step 4 of the $A_n$
algorithm above. In the dimension $4$ case
we again split the eigenspace according to the actions of the unipotent elements
already found.

However, in general in characteristic $2$, not every split toral subalgebra is
contained in a split maximal toral subalgebra, and therefore a heuristic algorithm
such as we suggest may fail to produce a Cartan subalgebra. It may produce a
maximal split toral subalgebra which is contained only in a non-split maximal
toral subalgebra.

\end{document}